\documentclass{lmcs} %%% last changed 2014-08-20
\pdfoutput=1

\usepackage{lastpage}
\lmcsdoi{18}{2}{12}
\lmcsheading{}{\pageref{LastPage}}{}{}%
{Mar.~29,~2021}{May~26,~2022}{}


\usepackage{hyperref}
\usepackage[utf8]{inputenc}
\newcommand{\N}{{\mathbb N}}
\newcommand{\Z}{{\mathbb Z}}

\DeclareMathOperator{\End}{End}
\DeclareMathOperator{\Pol}{Pol}

\newcommand{\alg}[1]{\mathbf{#1}}

\newcommand{\algA}{\alg{A}}
\newcommand{\algB}{\alg{B}}

\newcommand{\algD}{\alg{D}}
\newcommand{\algL}{\alg{L}}

\newcommand{\algU}{\alg{U}}

\DeclareMathOperator{\CSAT}{\mathsf{CSAT}}
\DeclareMathOperator{\CEQV}{\mathsf{CEQV}}

\DeclareMathOperator{\comP}{\mathsf{P}}
\DeclareMathOperator{\comNP}{\mathsf{NP}}

\DeclareMathOperator{\comcoNP}{\mathsf{coNP}}

\DeclareMathOperator{\CC}{CC}
\DeclareMathOperator{\ACC}{ACC}
\DeclareMathOperator{\AC}{AC}

\DeclareMathOperator{\MOD}{MOD}
\DeclareMathOperator{\AND}{AND}
\DeclareMathOperator{\OR}{OR}

\DeclareMathOperator{\arity}{ar}

\begin{document}

\title{CC-Circuits and the Expressive Power\texorpdfstring{\\}{ }of Nilpotent Algebras}
%\titlecomment{{\lsuper*}OPTIONAL comment concerning the title, \eg,
%  if a variant or an extended abstract of the paper has appeared elsewhere.}
\thanks{This work was supported by grant 18-20123S of the Czech Science Foundation, grant UNCE/SCI/022 of the Charles University Research Centre, and INTER-EXCELLENCE project LTAUSA19070 of the Czech Ministry of Education M\v{S}MT}	%optional

% affiliations are numbered automatically with a, b, c (see below)
% use the optional argument to indicate the affiliation(s) of each author
% omit the argument if there is only one author, or only one affiliation
\author[M.~Kompatscher]{Michael Kompatscher\lmcsorcid{0000-0002-0163-6604}}
%\author[B.~Name2]{Bob Name2}[a,b]
%\author[J.~Name3]{Josiah S.~Carberry\lmcsorcid{0000-0002-1825-0097}}[a]

% affiliation 1 (automatically numbered a)
\address{Department of Algebra, MFF, Charles University Prague, Czech Republic}	%optional
% write emails for all authors having that affiliation
\email{\texttt{michael@logic.at}}  %optional

%% etc.

%% required for running head on odd and even pages, use suitable
%% abbreviations in case of long titles and many authors:

%%%%%%%%%%%%%%%%%%%%%%%%%%%%%%%%%%%%%%%%%%%%%%%%%%%%%%%%%%%%%%%%%%%%%%%%%%%

%% the abstract has to PRECEDE the command \maketitle:
%% be sure not to issue the \maketitle command twice!

\begin{abstract}
  \noindent We show that CC-circuits of bounded depth have the same expressive power as circuits over finite nilpotent algebras from congruence modular varieties. We use this result to phrase a new algebraic version of Barrington, Straubing and Th\'erien's conjecture, which states that CC-circuits of bounded depth need exponential size to compute AND.

Furthermore we investigate the complexity of deciding identities and solving equations in a fixed nilpotent algebra. Under the assumption that the conjecture is true, we obtain quasipolynomial algorithms for both problems. On the other hand, if AND is computable by uniform CC-circuits of bounded depth and polynomial size, we can construct a nilpotent algebra in which checking identities is coNP-complete, and solving equations is NP-complete.
\end{abstract}

\maketitle

\section{Introduction}
Proving lower bounds on the size of Boolean circuits needed to compute explicit functions is a fundamental, but also notoriously hard problem in theoretical computer science. A rare exception with known sharp lower bounds is the parity function: In~\cite{furst-saxe-sipser} it was shown that parity cannot be computed by bounded depth circuits of polynomial size; \cite{yao-lowerbounds} then provided exponential lower bounds, which were further improved by H\r{a}stad~\cite{hastad-lowerbounds} to the almost optimal $2^{\Omega(n^{(d-1)^{-1}})}$, where $d$ stands for the depth of the circuits.

These results led to the question how much computational power we gain, if we also allow gates that describe parity or other counting functions in the construction of bounded depth circuits. More precisely, by such `counting gates' we mean $\MOD_m$-gates (for some $m \in \N$) of unbounded fan-in that output $1$, if the inputs sum up to $0$ modulo $m$, and $0$ otherwise (This use of $\MOD_m$ follows the notational convention used e.g. in~\cite{vollmer-circuitcomplexity, mckenzieNC1}. We remark that other authors such as~\cite{chattopadhyayETC, caussinus-depth2} actually denote the \emph{negation} of this function by $\MOD_m$). The class of functions that can be expressed by polynomially growing such bounded depth circuits is denoted by $\AC^0[m]$, and its union by $\ACC^0 = \bigcup_{m > 1} \AC^0[m]$.

An important step towards a characterization of $\ACC^0$ is to understand circuits that \emph{only} consist of $\MOD_m$-gates first. Such circuits are called $\CC[m]$-circuits. The functions that can be computed by bounded depth $\CC[m]$-circuits of polynomial size are denoted by $\CC^0[m]$, and their union over all $m$ by $\CC^0 = \bigcup_{m > 1} \CC^0[m]$. Despite being studied extensively, many questions about $\CC[m]$-circuits are still wide open. For instance their relationship for different values of $m$ is not well-understood, although this would be central for proving or disproving a conjecture of Smolensky~\cite{smolensky-conjecture}, which says that constant depth circuits having only $\AND$, $\OR$ and $\MOD_m$ gates cannot compute $\MOD_q$ in sub-exponential size, when $m$ is a prime power, and $q$ is coprime to it.

Another big open question is whether bounded depth $\CC[m]$-circuits can efficiently compute $\AND$ or not. If it is not the case, it would signify that there is a fundamental difference between circuits that use logical gates versus circuits that use counting gates. Barrington, Straubing and Th\'erien first conjectured that $\AND$ has exponential lower bounds, mirroring H\r{a}stad's result:

\begin{conj}[Barrington, Straubing and Th\'erien~\cite{barrington-NUDFA}] \label{conj:CCcircuits} 
Let $m \in \N$ and let $(C_n)_{n \in \N}$ be a family of bounded depth $\CC[m]$-circuits that compute $\AND$. Then $C_n$ is of size $2^{\Omega(n^q)}$ for some $q > 0$.
\end{conj}

For values of $m$ that have more than one prime divisor one can construct $\CC[m]$-circuits $(C_n)_{n \in \N}$ of depth $d$ and size $\mathcal O(2^{n^{(d-1)^{-1}}})$ that compute $\AND$ (see for instance~\cite{idziak-intermediate}). Thus one might be tempted to conjecture (in analogy to H\r{a}stad's result) that a value of $q$ close to $(d-1)^{-1}$ gives us optimal lower bound.  However, this is not the case: By the recent results in~\cite{CW-lowerbounds}, for every $q > 0$, there is a choice of $m$, such that there are depth $3$ $\CC[m]$-circuits of size  $2^{\mathcal O(n^q)}$ computing $\AND$. It was further shown in~\cite{IKK-modularcircuits} that for the correct choice of $m$ there are even $\CC[m]$-circuits of depth $2$ and size $2^{\mathcal O(n^q \log(n))}$ computing $\AND$. Thus any optimal bound in Conjecture~\ref{conj:CCcircuits} needs to depend not only on the depth $d$ but also on the modulus $m$.

Also weaker versions of Conjecture~\ref{conj:CCcircuits} can be found in the literature, e.g.~\cite{mckenzieNC1} conjectured that $\AND$ is not in $\CC^0$. Both the strong and the weak version of Conjecture~\ref{conj:CCcircuits} are open until today, with the best known general lower bound for $\AND$ being superlinear~\cite{chattopadhyayETC}. %NEXP result
It was further shown in~\cite{hansen-probabilisticCC0} that $\AND$ is in uniform \emph{probabilistic} $\CC^0$ (i.e. for every $c>0$, $\AND_n$ can be computed up to error $< \frac{1}{n^c}$ by $\comP$-uniform, bounded depth $\CC$-circuits with $\mathcal O(\log(n))$ random bits), which could be interpreted as evidence contrary to the conjecture.

However, for some special cases Conjecture~\ref{conj:CCcircuits} is confirmed: For prime powers $p^k$ it is well known that $\CC[p^k]$-circuits of bounded depth \emph{cannot} compute $\AND$~\cite{barrington-NUDFA}. The conjecture is further confirmed for the special case of bounded depth circuits with $\MOD_m$-gates at the input level, and $\MOD_{p^k}$-gates everywhere else, where $m \in \N$ and $p^k$ is a prime power (see~\cite{barrington-NUDFA}, or~\cite[Theorem~1]{caussinus-depth2}).

The first results about $\CC$-circuits arose from a characterization of them in the language of groups/monoids: In~\cite{barrington-NUDFA} Barrington, Straubing and Th\'erien used the notion of NUDFA (non uniform deterministic finite automata, introduced in~\cite{barrington-NUDFA-old}) to show that a language is in $\ACC^0$ if and only if it is accepted by a NUDFA over a solvable monoid, and in $\CC^0$ if and only if it is accepted by a NUDFA over a solvable group.\footnote{We remark that this characterization of $\CC^0$ was not explicitly stated in~\cite{barrington-NUDFA}, as the paper does not use the notion of $\CC$-circuits. However the relevant literature often attributes it to~\cite{barrington-NUDFA}, as it is a direct consequence of the presented results. One of the first explicit mentions of it can be found in~\cite[Theorem~2.8]{mckenzieNC1}. Also Conjecture~\ref{conj:CCcircuits} was originally stated as a conjecture about NUDFAs in~\cite{barrington-NUDFA}.} NUDFAs proved not only to be a fruitful tool in circuit complexity, but also led to new developments in algebra, regarding the study of equations in monoids and groups~\cite{barrington-eq-monoids}.

In this paper we give a new algebraic description of $\CC$-circuits, using concepts from universal algebra, more specifically \emph{commutator theory}. We show that, in some sense, $\CC$-circuits of bounded depth can be represented by the circuits over a nilpotent algebra from a congruence modular variety and vice-versa. As a corollary we obtain the following theorem:

\begin{thm} \label{theorem:main1}
Conjecture~\ref{conj:CCcircuits} is true, if and only if for every finite nilpotent algebra $\algA$ from a congruence modular variety, and for every $0 \in A$ there is a $q>0$, such that every sequence $(p_n)_{n \in \N}$ of nonconstant $0$-absorbing circuits $p_n(x_1,\ldots,x_n)$ over $\algA$ is of size $2^{\Omega(n^q)}$.
\end{thm}
Here, for a set $A$ and an element $0 \in A$ we call an operation $f\colon A^n \to A$ $0$-\emph{absorbing} if $0 = f(0,a_2,\ldots,a_n) = f(a_1,0,a_3,\ldots,a_n) = \cdots = f(a_1,\ldots,a_{n-1},0)$ holds for all values $a_1,\ldots,a_n \in A$.  Note that on the Boolean domain $A = \{0,1\}$ the only $n$-ary $0$-absorbing operations are the $n$-ary conjunction and the constant $0$-function. Thus nonconstant $0$-absorbing operations can be seen as a a natural generalization of conjunctions to arbitrary domains.  We remark that $0$-absorbing operations are of independent interest in commutator theory, as they characterize properties of the so called higher commutator~\cite{AichingerMudrinski}.

In the course of the proof of Theorem~\ref{theorem:main1} we are actually going to prove a stronger, but more technical result in Theorem~\ref{theorem:main}, which allows us to compute explicit bounds for $\AND$ from bounds on nonconstant $0$-absorbing circuits in nilpotent algebras, and vice-versa. We are further going to discuss how known results about $\CC[m]$-circuits correspond to known results about nilpotent algebras: the fact that for primes $p$, $\CC[p]$-circuits of bounded depth cannot compute $\AND$, corresponds for instance to the result that finite nilpotent algebras of prime power size have only nonconstant $0$-absorbing polynomial operations up to some arity~\cite{berman-blok}. The fact that Conjecture~\ref{conj:CCcircuits} holds for $\MOD_p$-$\MOD_q$ circuits~\cite{barrington-NUDFA} was essentially reproven in the language of nilpotent algebras in~\cite{IdziakKrzaczkowskiKawalek}.

At last we discuss the impact of Conjecture~\ref{conj:CCcircuits} on two computational problems, namely the circuit satisfaction problem $\CSAT(\algA)$ and the circuit equivalence problem $\CEQV(\algA)$ for fixed nilpotent algebras $\algA$. Here $\CSAT(\algA)$ models the decision problem, whether an equation over the algebra $\algA$ (encoded by two circuits) has a solution, while $\CEQV(\algA)$ asks, whether two circuits are equivalent. In~\cite{IdziakKrzaczkowski} Idziak and Krzaczkowski gave an almost complete complexity classification of both problems for algebras from congruence modular varieties, relating the complexity to the commutator theoretical properties of the given algebra. Essentially the only case left open are nilpotent, but not supernilpotent algebras (Problem 2 in~\cite{IdziakKrzaczkowski}). We show that, under the assumption of Conjecture~\ref{conj:CCcircuits}, there are quasipolynomial algorithms for both $\CSAT$ and $\CEQV$ of such algebras. On the other hand we show that, if $\AND$ is in \emph{$\comP$-uniform} $\CC^0$, then there is a nilpotent algebra with $\comNP$-complete circuit satisfaction problem, and $\comcoNP$-complete circuit equivalence problem.

We remark that, since this paper first appeared as a preprint, there were some new developments in the study of $\CSAT$ and $\CEQV$ for nilpotent algebras. Under the assumption of the Exponential Time Hypothesis (ETH), there are nilpotent algebras, for which both problems have proper quasipolynomial lower bounds~\cite{idziak-intermediate}. Together with the results in this paper, this implies that, under the ETH and Conjecture~\ref{conj:CCcircuits}, there are some nilpotent algebras, in which $\CSAT$ and $\CEQV$ can be solved in quasipolynomial, but not polynomial time. In particular no $\comP$/$\comNP$-c dichotomy result holds for $\CSAT$. This is an interesting contrast to the problem of solving \emph{systems} of equations in a given algebra $\algA$, since for every $\algA$, this problem is equivalent to a finite constraint satisfaction problem and hence either in $\comP$ or $\comNP$-complete by the CSP dichotomy theorem (independently proven by Zhuk~\cite{zhuk-dichotomy-short,zhuk-dichotomy} and Bulatov~\cite{bulatov-dichotomy}).  

The structure of this paper is as follows: In Section~\ref{sect:background} we discuss some necessary background from universal algebra and give a characterization of nilpotent algebras. In Section~\ref{sect:main} we prove Theorem~\ref{theorem:main1} and discuss its implications on $\CC$-circuits and nilpotent algebras. In Section~\ref{sect:CSAT} we derive the complexity results for $\CSAT$ and $\CEQV$ of nilpotent algebras.

\section{Background from universal algebra} \label{sect:background}

\subsection{Polynomials and circuits over algebras}
An algebra $\algA$ is a pair $(A,(f_i^{\algA})_{i \in I})$, where $A$ is a set (the \emph{universe} of $\algA$), and every element of the family $(f_i^{\algA})_{i \in I}$ is a finitary operation $f_i^{\algA} \colon A^{k_i} \to A$ (the \emph{basic operations} of $\algA$). The \emph{signature} of $\algA$ is the family $(k_i)_{i \in I}$, i.e. it describes which basic operations are of which arity. In this paper we are only going to consider \emph{finite algebras}, i.e. algebras that have a finite universe and only finitely many basic operations. By $\arity(\algA)$ we denote the maximal arity of the basic operations of $\algA$.

A \emph{term operation} of $\algA$ is an operation that can be obtained as a composition of basic operations of $\algA$. A \emph{polynomial operation} allows also the use of elements of $A$ in its construction. For the ring of integers $(\Z,+,0,-,\cdot)$ for instance, the polynomial operations are just the polynomial operations in the conventional sense (e.g.  the operation $p \colon \Z^3\to \Z$ defined by $p(x_1,x_2,x_3) = (2 x_1 \cdot x_3 + 4) \cdot x_1 \cdot x_2+1$). The set of all polynomial operations of $\algA$ is denoted by $\Pol(\algA)$.  Let $\algA$ and $\algB$ be two algebras with the same universe; if $\Pol(\algA) = \Pol(\algB)$ we say that $\algA$ and $\algB$ are \emph{polynomially equivalent}, if $\Pol(\algA) \subseteq \Pol(\algB)$, we say $\algB$ is a \emph{polynomial extension} of $\algA$.

Given a finite algebra, there are different ways of encoding its polynomial operations. The naive way is to just encode them by the string defining them. Such a string is usually referred to as a \emph{polynomial} over $\algA$. However, in an effort to compress the input, one can also consider \emph{circuits} over $\algA$, i.e. $A$-valued circuits with a unique output gate, whose gates are labelled by the basic operation of $\algA$.  We say that two circuits (respectively two polynomials) are \emph{equivalent}, if they define the same operation. 

By the following folklore result, the circuit approach does not only allow for a more concise, but also more stable representation of polynomial operations:

\begin{lem} \label{lemma:polynomialembedding}
Let $\algA$ and $\algB$ be two finite algebras with the same universe and $\Pol(\algA) \subseteq \Pol(\algB)$. Then we can rewrite every circuit $C(\bar x)$ over $\algA$ to an equivalent circuit $C'(\bar x)$ over $\algB$ in linear time.
\end{lem}

The circuit $C'$ in Lemma~\ref{lemma:polynomialembedding} can be easily obtained, by substituting every gate in $C$ by a circuit over $\algB$ defining the corresponding basic operation of $\algA$.  Note that it follows from Lemma~\ref{lemma:polynomialembedding} that the size of $C'$ is linear in the size of $C$. An analogue statement is however not true for the encoding by polynomials; in fact, the results in~\cite{HorvathSzabo-A4} imply that under the assumption $\comP\neq \comNP$ it even fails for quite well-behaved algebras such as the alternating group $A_4$. In this paper we are therefore only going to discuss the circuit encoding of polynomial operations. % (explaining the names CSAT and CEQV for our computational problems). %\footnote{However we remark (without giving a proof) that for the algebras used in Lemma~\ref{lemma:CCtonilpotent} and Proposition~\ref{proposition:main}, circuits can be rewritten to polynomials in polynomial time, and vice versa.Therefore also analogue statements to our results also hold for polynomials.}

For a circuit $C$ over $\algA$ with input gates $\bar x = (x_1,\ldots,x_n)$ we will write $C(\bar x)$ both for the circuit, and the $n$-ary polynomial operation induced by it, but this should never cause any confusion. We call a circuit $C$ constant/$0$-absorbing/etc., if the polynomial operation $C(\bar x)$ defined by it is constant/$0$-absorbing/etc.

\subsection{Nilpotent algebras from congruence modular varieties}
Commutator theory is a field of universal algebra that generalizes concepts from the commutator theory of groups to arbitrary algebras. In particular there is the notion of a central series of congruences, which allows to define nilpotent algebras (as a generalization of nilpotent groups and nilpotent rings). Since we are only interested in nilpotent algebras from congruence modular varieties, we refrain from giving the original definition (via the so called \emph{term congruence commutator}) and refer to~\cite{FreeseMcKenzie-Commutator} for background.

Algebras from \emph{congruence modular varieties} form a quite broad class that contains many examples of interest in abstract algebra and computer science, such as: lattices, Boolean algebras, fields, rings, modules, groups, quasigroups, and all their polynomial extensions. Commutator theory works especially well in the congruence modular case.

In a congruence modular variety, an algebra is \emph{Abelian} (or \emph{1-nilpotent}) if and only if it is polynomially equivalent to a \emph{module}. Here, a module $\algA$ is an algebra $\algA = (A,+,0,-, (r)_{r \in R})$, where $(A,+,0,-)$ is the underlying Abelian group and every ring element $r \in R$ is considered as a unary operation $r(x) = r \cdot x$. Therefore all operations of such an Abelian algebra $\algA$ are affine operations $p(x_1,\ldots,x_n) = \sum_{i = 1}^n r_i x_i + c$ for some scalars $r_i \in R$ and a constant $c \in A$.

General nilpotent algebras from congruence modular varieties can, in some sense, be decomposed into Abelian algebras (this follows from Proposition~7.1. in~\cite{FreeseMcKenzie-Commutator}). We are going to treat this characterization as a definition of nilpotence:

\begin{defi} \label{theorem:freese-mckenzie}
An algebra $\algA = (A,(f^{\algA})_{f \in I})$ from a congruence modular variety is
\begin{itemize}
\item \emph{1-nilpotent} (or \emph{Abelian}) if and only if it is polynomially equivalent to a module;  
\item \emph{$n$-nilpotent}, if there are algebras $\algL = (L,(f^{\algL})_{f \in I})$ and $\algU = (U,(f^{\algU})_{f \in I})$ in the same signature as $\algA$ such that 
\begin{itemize}
\item $\algL$ is Abelian and $\algU$ is $(n-1)$-nilpotent,
\item $A = L \times U$, where $L$ and $U$ are the universes of $\algL$ and $\algU$ respectively,
\item Every basic operation $f^{\algA}$ of $\algA$ is of the form\\ 
$f^{\algA}((l_1,u_1),\ldots,(l_k,u_k)) = (f^\algL(l_1,\ldots,l_k) + \hat f(u_1,\ldots,u_k), f^\algU(u_1,\ldots,u_k)),$\\ for a function $\hat f \colon U^k \to L$.
\end{itemize}
\end{itemize}
Here $+$ denotes the addition of the module that is polynomially equivalent to $\algL$. We write $\algA = \algL \otimes \algU$ for the decomposition of $\algA$ into $\algL$ and $\algU$.
\end{defi}

When talking about nilpotent algebras in this paper, we will from now on always implicitly assume that they are from a congruence modular variety. 

Note that by Definition~\ref{theorem:freese-mckenzie}, $n$-nilpotent algebras can be regarded as an $(n-1)$-nilpotent algebra $\algU$ `acting' on the Abelian algebra $\algL$ by the operations $\hat f \colon U^k \to L$. The \emph{degree of nilpotency} of $\algA$ is the smallest $n$, such that $\algA$ is $n$-nilpotent. Note that this degree is always bounded by $\lfloor \log_2|A| \rfloor$.
 %Note that, on a  conceptual level, this reflects the structure of $\CC$-circuits of bounded depth, or also the wreath product construction that was used in~\cite{barrington-NUDFA}. 

By a recent result of Aichinger every nilpotent algebra has a nicely behaved extension by some Abelian group operations:

\begin{thm}[Corollary of Theorem~4.2. in~\cite{aichinger-spectrum}] \label{theorem:coordinatisation}
Let $\algA$ be a finite nilpotent algebra from a congruence modular variety and let $0 \in A$. Then there exists a nilpotent algebra $\algB$ with the same universe $A$, such that 
\begin{itemize}
\item $\Pol(\algA) \subseteq \Pol(\algB)$,
\item $\algB = \algL_1 \otimes (\algL_2 \otimes ( \cdots  \otimes \algL_n ))$, such that 
\begin{itemize}
\item every $\algL_i$ is polynomially equivalent to a module with group $(L_i,+,0,-) = \Z_{q_i}^{k_i}$ and a ring $R \leq \End(\Z_{q_i}^{k_i})$, for a prime $q_i$ and $k_i \in \N$,
\item $\algB$ contains basic operations $+$, $0$, $-$ such that $(A,+,0,-) = \prod_{i = 1}^n (L_i,+,0,-)$.
\end{itemize}
\end{itemize}
\end{thm}

\noindent We remark that Theorem~4.2. in~\cite{aichinger-spectrum} is phrased in a different way from our statement: it states that to every central series of congruences $0_A = \alpha_0 < \alpha_1 < \ldots < \alpha_n = 1_A $ of $\algA$ and every constant $0 \in A$, we can associate an Abelian group $(A,+,0,-)$,  such that $(\alpha_i)_{i=0}^n$ is still a central series of $\algB$, the extension of $\algA$ by $+,0,-$. By applying Proposition~7.1.  from~\cite{FreeseMcKenzie-Commutator} iteratively to $\algB$, we then obtain the algebras $\algL_i$ such that $\algB= \algL_1 \otimes (\algB / \alpha_1) =  \algL_1 \otimes (\algL_2 \otimes (\algB / \alpha_2)) = \cdots = \algL_1 \otimes (\algL_2 \otimes ( \cdots  \otimes \algL_n))$ such that $(A,+,0,-) = \prod_{i = 1}^n (L_i,+,0,-)$.  If we pick the central series to be of maximal length,  it holds that $(L_i,+,0,-) = \Z_{q_i}^{k_i}$ for some prime $q_i$ by Lemma 4.1. in~\cite{aichinger-spectrum}. 

By Lemma~\ref{lemma:polynomialembedding} every circuit over $\algA$ can be rewritten in linear time to a circuit over its polynomial extension $\algB$ given by Theorem~\ref{theorem:coordinatisation}; working in the extension $\algB$ instead of $\algA$ will simplify our proof.  In particular, we are going to use that every element $x \in A$ can be identified with a tuple in $(A,+,0,-) = \prod_{i = 1}^n (L_i,+,0,-) = \prod_{j=1}^s  \Z_{p_j}$ such that $|A| = \prod_{j=1}^s {p_j}$ is the prime decomposition of $|A|$ (possibly containing repetitions). 

We want to stress here, that only the group $(A,+,0,-)$ can be decomposed into a direct product of prime sized algebras,  but not the entire algebra $\algB$.  Moreover, the Abelian algebras $\algL_i$ can in general not be further decomposed with respect to $\otimes$; an example is the Abelian algebra $(\Z_2 \times \Z_2, +,0,-, (r)_{r \in R})$, where $R = \End(\Z_2 \times \Z_2)$ is the full endomorphism ring.

\section{The equivalence of CC-circuits and circuits over nilpotent algebras} \label{sect:main}

Our proof of Theorem~\ref{theorem:main1} is based on expressing circuits over a fixed nilpotent algebra $\algA$ as $\CC[m]$-circuits, and vice-versa. It is however a priori not clear if and how this is possible, as $\CC[m]$-circuits are Boolean valued, whereas the universe of $\algA$ can be arbitrary. Because of this, most of the time we are not going to work with $\CC[m]$-circuits themselves, but an $m$-valued analogue, which we call $\CC^+[m]$-circuits. We define $\CC^+[m]$-circuits and discuss some of their properties in the next subsection. This is then followed by the proof and discussion of Theorem~\ref{theorem:main}, which in particular implies Theorem~\ref{theorem:main1}.

\subsection{$\CC^+[m]$-circuits} 

\begin{defi}
A \emph{$\CC^+[m]$-circuit} $C$ is a circuit over the set $\Z_m = \{0,1,\ldots,m-1\}$ containing $+$-gates, constant gates,  and $\MOD_m$-gates of arbitrary fan in.  We  interpret $+$ as addition modulo $m$.  As in the Boolean case, $\MOD_m$-gates output $1$, if their inputs sum up to $0$ modulo $m$, and $0$ otherwise.
\end{defi}

Without loss of generality we can assume that in such a $\CC^+[m]$-circuit $C$ there are at most $m$ wires between two gates (otherwise we take their number modulo $m$). Furthermore, every wire from a $+$-gate to some $\MOD_m$-gate can be substituted by wires that go directly from the inputs of the $+$-gate to the $\MOD_m$-gate.  By doing this substitution for all such wires (and reducing their number in each step modulo $m$) we obtain an equivalent circuit, whose size is linear in the size of the original. Thus, from now on, we assume that every $\CC^+[m]$-circuit has at most one $+$-gate, and this is the output gate.

We next discuss the relationship between nonconstant $0$-absorbing $\CC^+[m]$- and $\CC[m]$-circuits. Note that an $n$-ary $\CC[m]$-circuit $C$ is nonconstant and $0$-absorbing if and only if $C(\bar x) = 1$ if $\bar x = (1,1,\ldots,1)$ and $C(\bar x) = 0$ else.  So it defines the $n$-ary $\AND$ operation.

\begin{lem} \label{lemma:CC+} \
\begin{enumerate}
%\item Every $\CC^+[m]$-circuit $C$ of depth $d$ can be rewritten in linear time to a $\CC^+[m]$-circuit $C'$ of depth $\leq d$, which has no $+$-gates except at its output.
\item For every nonconstant $0$-absorbing $\CC^+[m]$-circuit $C(\bar x)$ of depth $d$ we can compute in linear time a $\CC[m]$-circuit $C'(\bar x)$ of depth $d$ that defines $\AND$.
\item For every $\CC[m]$-circuit $C(\bar x)$ of depth $d$ defining $\AND$ we can compute in linear time a nonconstant $0$-absorbing $\CC^+[m]$-circuit $C'(\bar x)$ of depth $d+1$.
\end{enumerate}
\end{lem}

\begin{proof} %(1) follows straightforward from the fact that a wire from a $+$-gate to some $\MOD_m$-gate can be substitutes by wires from the inputs of the $+$-gate to the $\MOD_m$-gate.\\
For (1), let $C(x_1,\ldots,x_n)$ be a $\CC^+[m]$-circuit of depth $d$ that defines a nonconstant $0$-absorbing function. Thus there is some tuple $\bar a \in (\Z_m \setminus \{0\})^n$ such that $C(\bar a) \neq 0$. Without loss of generality we can assume that $\bar a = (1,1,\ldots,1)$, otherwise we duplicate every wire connected to the input gate $x_i$ $a_i$-many times. Without loss of generality we can further assume that $C$ has only constant gates labelled by $1$ (otherwise we substitute the constant gate $c$ by $c$-many copies of $1$). If $C$ contains no $+$-gate at all, we set $C = C'$ and are done. If $C$ has a $+$-gate as an output gate that sums up  $g_1,\ldots,g_k$ we construct $C'$ by substituting this $+$-gate by a $\MOD_m$-gate that has $g_1,\ldots,g_k$ and $m - C(\bar a)$ as input.

For (2) note that $C$ itself might not be $0$-absorbing when evaluated over $\Z_m$. However the circuit $C'(\bar x) = C(\MOD_m(1-x_1),\ldots,\MOD_m(1-x_n))$ is.
\end{proof}

As a consequence of Lemma~\ref{lemma:CC+} (1), lower bounds on the size of $\CC[m]$-circuit of depth $d$ defining $\AND$ are also lower bounds on the size of nonconstant $0$-absorbing $\CC^+[m]$-circuit of depth $d$. By Lemma~\ref{lemma:CC+} (2) also the reverse statement holds, up to decreasing the depth by 1. 

We continue by discussing representations of arbitrary functions by $\CC^+[m]$-circuits. A useful concept for it is the \emph{essential depth} of a $\CC^+[m]$-circuit $C$, which we define to be the depth of $C$, not counting its $+$-gates.

\begin{lem} \label{lemma:CC+fct} \
\begin{enumerate}
\item If $m>2$, then every function $f \colon \Z_m^n \to \Z_m$ can be represented by a $\CC^+[m]$-circuit of essential depth $\lceil \log_2(n) \rceil + 1$.
\item If $m$ has two distinct prime factors, then every function $f\colon \Z_m^n \to \Z_m$ can be represented by a $\CC^+[m]$-circuit of essential depth $3$. 
\end{enumerate}
\end{lem}

\begin{proof}
In order to prove (1), let us define the series of circuits 
\begin{align*}
C_1(x_1) &= \MOD_m(x_1), \text{ and }\\
C_{2^{k+1}}(x_1,\ldots,x_{2^{k+1}}) &= \MOD_m(C_{2^k}(x_1,\ldots,x_{2^k}), C_{2^k}(x_{2^k+1},\ldots,x_{2^{k+1}}), m - 2),
\end{align*}
for every $k \in \N$. Note that $C_{2^{k+1}}$ outputs $1$ if all inputs are equal to $0$ and $0$ otherwise.  We can further construct circuits $C_n$ with this behaviour for all arities $2^k < n \leq 2^{k+1}$ by identifying some of the variables of $C_{2^{k+1}}$. Now (1) follows from the fact that every function $f \colon \Z_m^n \to \Z_m$ can be obtained as the sum of translations of $C_n$ by constant tuples, i.e. $f(\bar x) = \sum_{\bar a \in A^n} f(\bar a) \cdot C_n(\bar x - \bar a)$ (where $f(\bar a) \cdot y$ is a shorthand for the sum of $f(\bar a)$-many copies of $y$).

For (2) recall that all Boolean operations, in particular the $n$-ary $\AND$, can be written as $\CC[m]$-circuits of depth 2 (see~\cite{barrington-NUDFA}). Thus $C_n(\bar x) = \AND (\MOD_m(x_1),\ldots,\MOD_m(x_n))$ is a $\CC^+[m]$-circuit of depth $3$ that describes the characteristic function of $(0,\ldots,0)$. As in the proof of (1), we can obtain all other functions as linear combinations $f(\bar x) = \sum_{\bar a \in A^n} f(\bar a) \cdot C_n(\bar x - \bar a)$.
\end{proof}

\subsection{The main result}
Recall that in Definition~\ref{theorem:freese-mckenzie} we considered the operations of $n$-nilpotent algebras as `actions' of an $(n-1)$-nilpotent algebra on an Abelian one. Also in a $\CC^+[m]$-circuit we can think of a $\MOD_m$-gates as receiving inputs from gates of lower depth and having an output in the Abelian group $\Z_m$. This viewpoint allows us  to construct a nilpotent algebra in which we can interpret all $\CC^+[m]$-circuits of fixed depth $d$:

\begin{lem} \label{lemma:CCtonilpotent}
For all $m,d \in \N$ there is a $(d+1)$-nilpotent algebra $\algD = (\Z_m^{d+1},+,0,-, h, *, (f_i)_{i=1}^d)$, such that for every $\CC^+[m]$-circuit $C$ of essential depth $d$, there is a circuit $C'$ over $\algD$ with $C'(x_1,\ldots,x_n) = (C(\pi_{d+1}(x_1),\ldots,\pi_{d+1}(x_n)),0,\ldots,0)$ (where $\pi_i$ denotes the projection of $D$ to the $i$-th component of $(\Z_m)^{d+1}$). Further $C'$ can be computed from $C$ in linear time.
\end{lem}

\begin{proof}
As usual, $+,0,-$ denote the Abelian group operation on $\Z_m^{d+1}$.  Further we define the unary operation $h(x) = (0,\ldots,0,\pi_{d+1}(x))$, the binary operation $x_1 *x_2 = (\sum_{i=1}^{d+1}\pi_{i}(x_1) + \pi_{i}(x_2),0,\ldots,0)$, and the unary operations $(f_i)_{i = 1}^d$ by 
$$f_{i}(x) = \begin{cases} (0,\ldots,0,1,0,\ldots,0) \text{ if } \pi_{i+1}(x) + \pi_{i+2}(x) + \cdots + \pi_{d+1}(x) = 0,\\(0,0,\ldots,0) \text{ else;} \end{cases}$$
here $1$ lies on the $i$-th coordinate.
Note that $\algD = \algL_1 \otimes (\algL_2 \otimes \cdots (\algL_{d} \otimes \algL_{d+1}) \cdots)$, such that every $\algL_i$ is polynomially equivalent to $\Z_m$.  Since $\algD$ contains group operations, and groups have modular congruence lattices, $\algD$ generates a congruence modular variety. Thus, by Definition~\ref{theorem:freese-mckenzie}, $\algD$ is $d+1$-nilpotent. Let $C(y_1,\ldots,y_n)$ now be a $\CC^+[m]$-circuit of essential depth $d$. Without loss of generality we can assume that $C(y_1,\ldots,y_n)$ has at most one $+$-gate at its output level. 

By the depth $r$ of a gate $g$ of $C$ we denote the distance of $g$ from one of the input or constant gates.  Further let $C_g$ be the subcircuit of $C$, consisting of all vertices that have a directed path to $g$.  We are going to show by induction over the depth $r = 0,1,\ldots,d$ that for every gate $g$ of depth $r$ there is a circuit $C_g'$ over $\algA$, such that $\pi_{d+1-r} C_g'(x_1,\ldots,x_n) = C_g(\pi_{d+1}(x_1),\ldots,\pi_{d+1}(x_n))$ and $\pi_{j} C_g'(x_1,\ldots,x_n)= 0$,  for all $j \neq d+1-r$.  Note that for $r=d$, this is essentially the statement of the lemma.

The gates $g$ of depth $r=0$ are either input gates or constant gates.  For every constant gate $C_g =c$ of $C$ we set $C_g' = (0,0,\ldots,c)$  and every input gate $C_g = y_i$ of $C$ we set $C_g' = h(x_i)$.  Clearly $C_g'$ then satisfies the our claim.

For an induction step $r-1 \to r$, let $g$ be a gate of depth $r$.  Note that $g$ must be a $\MOD_m$-gate, and thus $C_g = \MOD_m(C_{g_1}+ C_{g_2} + \ldots +C_{g_n})$, where all $g_i$ are gates of depth $r-1$ or lower.  By induction hypothesis, for every $g_i$ there exists a circuit $C_{g_i}'$ over $\algA$, whose projection to the component of $\Z_m^{d+1}$ corresponding to the depth of $g_i$ is $C_{g_i}(\pi_{d+1}(x_1),\ldots,\pi_{d+1}(x_n))$, and whose projection to every other component is $0$.We then define $C_g' = f_{d+1-r}(C_{g_1'} + C_{g_2'} + \ldots +C_{g_n'})$. By the definition of the function $f_{d+1-r}$,  this circuit satisfies $\pi_{d+1-r} C_g'(x_1,\ldots,x_n) = C_g(\pi_{d+1}(x_1),\ldots,\pi_{d+1}(x_n))$ and $\pi_{j} C_g' = 0$ for all $j \neq d+1-r$, which is what we wanted to prove. 

If $C$ has a $\MOD_m$-gate $g$ of depth $d$ as output gate, then we simply set $C' = C_g'$ and are done.  In the other case $C'$ has a $+$-gate (of depth $d+1$) as output gate.  Hence $C = C_{g_1} + C_{g_2} + \ldots + C_{g_n}$,  for gates $g_1,\ldots,g_n$ of depth at most $d$. Let $C' = C_{g_1}' *C_{g_2}' * \ldots *C_{g_n}'$.  By definition of $*$, this $C'$ also satisfies the statement of the lemma. % The textual modification in this paragraph allows to keep the definitions of circuit satisfiability and circuit equivalence on the same page (p.11).
\end{proof}

The other direction, i.e. showing that circuits over a fixed nilpotent algebra $\algA$ can be expressed by bounded depth $\CC^+[m]$-circuits, requires some more work.

Recall that by Theorem~\ref{theorem:coordinatisation} every nilpotent algebra $\algA$ can be extended by group operations $+,0,-$ such that $(A,+,0,-)$ is the direct product of prime order groups $\prod_{i = 1}^s \Z_{p_i}$, with $|A| = \prod_{j = 1}^s p_j$.  So we can identify an element $a \in A$ with the tuple $(\pi_1(a)\ldots,\pi_s(a)) \in \prod_{j = 1}^s \Z_{p_j}$.  Let $m$ be the product of all \emph{distinct} prime factors of $|A|$. Every factor $\Z_{p_j}$ clearly embeds into $\Z_m$ (by the map $x \mapsto (mp_j^{-1}) x$). Therefore there is also a natural embedding $e \colon (A,+,0,-) = \prod_{j = 1}^s \Z_{p_j} \to (\Z_m)^s$.

Our goal is to interpret circuits over the algebra $\algA$ as $\CC^+[m]$-circuits using this embedding. So for a circuit $C$ over $\algA$, we want to find a $\CC^+[m]$-circuit $C'$ such that $C'(e(x_1),\ldots,e(x_k)) = e C(x_1,\ldots,x_k)$.  Note that for every input gate $x_i$ of $C$, the $\CC^+[m]$-circuit $C'$ needs to have $s$ input gates for $e(x_i) = (\pi_1(x_i)\ldots,\pi_s(x_i))$.  For the same reason $C'$ will also have $s$ output gates.

\begin{prop} \label{proposition:main}
Let $\algA$ be a finite nilpotent algebra, and $\algB$ be its extension by the Abelian group operations $+,0,-$, given by Theorem~\ref{theorem:coordinatisation}. \\
Let $m > 2$ be the product of prime divisors of $|A|$, and $e \colon (A,+,0,-) \to (\Z_m)^s$ the natural group embedding. Then there is a $d \in \N$ such that for every circuit $C$ over $\algB$ there is a $\CC^+[m]$-circuit $C'$ of essential depth $d$ with $C'(e(x_1),\ldots,e(x_k)) = e C(x_1,\ldots,x_k)$.
\begin{enumerate}
\item In general $d \leq (\lceil \log_2 (s \cdot \arity(\algB)) \rceil +1 ) \cdot (\log_2|A|-1)$,
\item If $m$ is not a prime, then $d \leq 3(\log_2|A|-1)$.
\end{enumerate}
Furthermore $C'$ can be computed from $C$ in linear time.
\end{prop}

\begin{proof}
Recall that, by Theorem~\ref{theorem:coordinatisation},$$\algB = \algL_1 \otimes (\algL_2 \otimes \cdots (\algL_{n-1} \otimes \algL_n) \cdots),$$ such that every $\algL_i$ is polynomially equivalent to a module, whose group is $\Z_{p_i}^{k_i}$ for a prime $p_i$. With respect to this representation, every basic operation of $\algB$ is of the form 
$$f = (f^{\algL_1} + \hat f_1, f^{\algL_2} + \hat f_2, \ldots, f^{\algL_{n-1}} + \hat f_{n-1}, f^{\algL_{n}}),$$
where $f^{\algL_i}$ is an operation of $\algL_i$ and $\hat f_i$ only depends on the projection of $A$ to $L_{i+1} \times \cdots \times L_n$. Without loss of generality we can assume that each basic operation of $\algB$ is either `of affine type', meaning that $f = (f^{\algL_1}, f^{\algL_2}, \ldots, f^{\algL_{n-1}}, f^{\algL_{n}})$, or of `hat type' meaning that $f = (0,\ldots, 0,\hat f_j,0,\ldots,0)$.  If this is not the case we substitute $f$ by the basic operations $(f^{\algL_1}, \ldots, f^{\algL_{n}})$ and $(0,\ldots,0,\hat f_{j},0,\ldots,0)$ for every $j=1,\ldots,n-1$. The resulting algebra is also $n$-nilpotent and a polynomial extension of $\algB$.  Note further that this extension has the same maximal arity as $\algB$.

Every basic operation $f^{\algL_i} \colon L_i^l \to L_i$ of $\algL_i$ is equal to an affine combination $f^{\algL_i}(z_1,\ldots,z_l) = \sum_{j = 1}^l r_i \cdot z_i + c$, where the $r_i$ are linear maps. Therefore we can consider $f^{\algL_i}$ as an affine map from the vector space $\Z_{p_i}^{l \cdot k_i}$ to $\Z_{p_i}^{k_i}$. As such, it can be represented by a circuit over $\Z_{p_i}$ with $lk_i$ inputs and $k_i$ outputs, consisting only of $+$-gates. By duplicating every wire in this circuit $mp^{-1}$-many times, we obtain a corresponding $+$-circuit over $\Z_m$. By doing this for every projections $f^{\algL_i}$ of a basic operation $f =(f^{\algL_1}, \ldots, f^{\algL_{n}})$ of `affine type', we can show that there is a $\CC^+[m]$ circuits $C_f$ only consisting of $+$-gates, such that $C_f(e(x_1),\ldots,e(x_l)) = e f(x_1,\ldots,x_l)$.

Next, for every operation of `hat type' $f \colon A^l \to A$, we take an arbitrary extension of it to a function $g \colon (\Z_m^s)^l \to \Z_m^s$.  Since $l \leq \arity(\algB)$,  by Lemma~\ref{lemma:CC+fct} (1), this $g$ can be represented by a $\CC^+[m]$-circuit $C_f$ of essential depth $\lceil \log_2 (s\cdot \arity(\algB)) \rceil +1$. If $m$ is the product of two or more prime factors, then by Lemma~\ref{lemma:CC+fct} (2), $C_f$ can even be written as a $\CC^+[m]$-circuits of essential depth $3$.

Now let $C$ be an arbitrary circuit of the algebra $\algB$.  For every basic operation of hat-type $f$, there is a $j \in \{1,\ldots,n-1\}$ such that $\hat f_j$ only depends on $L_{j+1} \times \cdots \times L_n$,  and all other coordinates are $0$.  For any other operation $g$ that is of of hat-type with respect to a smaller index $p \leq j$ the identity $f(x_1,\ldots,x_k + c \cdot g(\bar y), \ldots,x_l) = f(x_1,\ldots,x_k, \ldots,x_l)$ holds, for all $k = 1,\ldots,l$.  Using these identities, the circuit $C$ can be rewritten in linear time to a circuit, in which no directed path contains more than $n-1$ gates of hat type.

If we then substitute all gates $f$ of $C$ by the corresponding $\CC^+[m]$-circuits $C_f$ as described above, we obtain a $\CC^+[m]$-circuit $C'$ of essential depth $(\lceil \log_2 (s\cdot \arity(\algB)) \rceil +1)(n-1)$. In the case, in which $m$ has two or more prime factors,  the same argument gives us a circuit of essential depth $3(n-1)$.

Since $n \leq \log_2|A|$, this proves the proposition.
\end{proof}

Note that the $sk$-ary $\CC^+[m]$-circuit $C'(y_1, \ldots,y_{sk})$ in Proposition~\ref{proposition:main} does not need to be $0$-absorbing, if $C$ is $0$-absorbing. However we can obtain a $0$-absorbing circuit with the same range as $e C$ by taking the circuit $C''(y_1, \ldots,y_{sk})$ that computes $C'(c_1 y_1, c_2 y_2, \ldots, c_s y_{sk})$, with $c_i = m p_i^{-1}$. 

We are now ready to prove the main result.

\begin{thm} \label{theorem:main}
For a fixed nilpotent algebra $\algA$ and $0 \in A$, let $f_\algA(n)$ denote the minimal size of a nonconstant $0$-absorbing $n$-ary circuit over $\algA$. For two integers $d,m$, let $g_{m,d}(n)$ be the smallest size of an $n$-ary $\CC[m]$-circuit of depth $d$ computing $\AND$. Then: \\

\begin{enumerate}
\item Let $\algA$ be a finite nilpotent algebra, and $m$ be the product of prime divisors of $|A|$. If $m$ is not a prime, then $f_\algA(n) \geq K g_{m,d}(n)$ for $d = 3 \lfloor \log_2|A| \rfloor - 2$ and some $K>0$.
\item Vice versa, for every $d > 1$ and every product of two or more distinct primes $m$, there is a $(d+2)$-nilpotent algebra $\algB$ with $g_{m,d}(n) \geq K' f_{\algB}(n)$ for some $K' > 0$.
\end{enumerate}
\end{thm}

\begin{proof}
To see (1), let $(C_n)_{n \in \N}$ be a sequence of nonconstant $0$-absorbing circuits over $\algA$ such that $|C_n| = f_\algA(n)$. We can regard every $C_n$ also as a circuit over the nilpotent extension $\algB$ of $\algA$ given by Theorem~\ref{theorem:coordinatisation}. This extension $\algB$ is nilpotent and contains a group operation $+$ such that $(A,+,0,-) = \prod_{j = 1}^s \Z_{p_j}$. By Proposition~\ref{proposition:main}, for every $n$ there is a $s n$-ary $\CC^+[m]$-circuit $C_n'$ of essential depth at most $3 \lfloor \log_2|A| \rfloor - 3$ such that $C_n'(e(x_1),\ldots,e(x_n)) = e C_n(x_1,\ldots,x_n)$. Note that, by eliminating all $+$-gates, except the output gates, we obtain a $\CC^+[m]$-circuit $C_n'$ of depth at most $d = 3 \lfloor \log_2|A| \rfloor - 2$.

Since $C_n$ is $0$-absorbing, the circuit $C_n'(c_1 y_1,\ldots,c_{s} y_{sn})$, for $c_i = (mp_i^{-1})$ is $0$-absorbing and has the same range as $e C_n$. As $C_n$ is not constant, there is an output gate of $C_n'$ that induces a nonconstant $0$-absorbing operation. By Lemma~\ref{lemma:CC+} (1) we can compute in linear time a $\CC[m]$-circuit of the same depth, that defines the $sn$-ary $\AND$. Therefore there is a $K$, such that $f_\algA(n) \geq K g_{m,d}(s n) \geq K g_{m,d}(n)$; the constant $K$ results from the fact that all computations only required linear time.

For (2) note that by Lemma~\ref{lemma:CC+} (2), for every $\CC[m]$-circuit $C$ of depth $d$ defining $\AND$, we can construct in linear time a nonconstant, $0$-absorbing $\CC^+[m]$-circuit $C'$ of depth $d+1$. Let us define $\algB$ to be the $(d+2)$-nilpotent algebra $\algD$ given by Lemma~\ref{lemma:CCtonilpotent} for depth $d+1$. We then can compute a circuit $C''$ over $\algB$ such that $C''(x_1,\ldots,x_n)$ evaluates to $(C'(\pi_{d+2}(x_1),\ldots,\pi_{d+2}(x_n)),0,\ldots,0)$. $C''$ is clearly nonconstant and $0$-absorbing. This concludes the proof of (2).
\end{proof}

Theorem~\ref{theorem:main} directly implies that Conjecture~\ref{conj:CCcircuits} is true, if and only if nonconstant $0$-absorbing polynomial operations in nilpotent algebras require circuits of exponential size. So we obtain Theorem~\ref{theorem:main1} as a direct corollary.

Moreover Theorem~\ref{theorem:main} implies that also the weak version of Conjecture~\ref{conj:CCcircuits}, which states that $\AND$ is not in $\CC^0$, has an algebraic counterpart by the following equivalence:

\begin{cor} \label{corollary:ANDinCC0}
$\AND$ is in $\CC^0$ if and only if there is a nilpotent algebra $\algA$ (from a congruence modular variety) and a series of circuits $(C_n(\bar x))_{n \in \N}$ over $\algA$ that grows polynomially and defines nonconstant $0$-absorbing polynomials of all arities.
\end{cor}

We further remark that, if $m$ is an odd prime, Proposition~\ref{proposition:main} allows us to reprove results about nilpotent algebras: Since we know that bounded depth $\CC[m]$-circuits are not able to define $\AND$, Proposition~\ref{proposition:main} implies that finite nilpotent algebras of prime power order only have nonconstant $0$-absorbing polynomials up to some arity. This was independently already shown in~\cite{berman-blok}. In fact, finite nilpotent algebras that have nonconstant $0$-absorbing polynomials only up to some fixed arity (so called \emph{supernilpotent} algebras) are characterized by being direct products of nilpotent algebras of prime power size~\cite{kearnes-spectra}. Nilpotent groups and rings are examples of supernilpotent algebras.

In~\cite{IdziakKrzaczkowskiKawalek} it was further shown that, if $\algA$ is 2-nilpotent and $\algA = \algL \otimes \algU$, where $\algL$ and $\algU$ are two vector spaces of different characteristics $p$ and $q$,  there are nonconstant $0$-absorbing polynomial operations that can be computed by circuits of exponential size $\mathcal O (c^n)$; however they cannot be computed by circuits of smaller size.  In the spirit of Theorem~\ref{theorem:main}, this can be also seen as the algebraic counterpart of the well-known fact that $\MOD_p$-$\MOD_q$ circuits need at least size $\Omega(c^n)$ to compute $\AND$~\cite{barrington-NUDFA}.

\section{Circuit satisfiability and equivalence} \label{sect:CSAT}
In this section we discuss the complexity of the circuit satisfiability and the circuit equivalence problem for nilpotent algebras. The circuit satisfiability problem $\CSAT(\algA)$ models the question, whether a single equation over the algebra $\algA$ has a solution; the circuit equivalence problem $\CEQV(\algA)$ asks whether an equation holds for all assignments of variables. Both problems were introduced in~\cite{IdziakKrzaczkowski} and are formally defined as follows:\\

\textsc{Circuit satisfiability} $\CSAT(\algA)$\\
\textsc{Input:} Two circuits $C,C'$ over $\algA$ with input gates $\bar x = (x_1,\ldots,x_n)$\\
\textsc{Question:} Is there a tuple $\bar a \in A^n$ such that $C(\bar a) = C'(\bar a)$?\\

\textsc{Circuit equivalence} $\CEQV(\algA)$\\
\textsc{Input:} Two circuits $C,C'$ over $\algA$ with input gates $\bar x = (x_1,\ldots,x_n)$\\
\textsc{Question:} Is $C(\bar a) = C'(\bar a)$ for all $\bar a \in A^n$?\\

Circuits over finite algebras can be evaluated in polynomial time. Therefore $\CSAT(\algA)$ is always in $\comNP$ and $\CEQV(\algA)$ in $\comcoNP$. The major question then is, for which algebras $\CSAT$ and $\CEQV$ are tractable, and for which they are $\comNP$-complete, respectively $\comcoNP$-complete. In particular this is still open for nilpotent algebras from congruence modular varieties. We first show that - under the assumption that Conjecture~\ref{conj:CCcircuits} is true - there are algorithms for both $\CEQV(\algA)$ and $\CSAT(\algA)$ that run in quasipolynomial time. This gives us a conditional answer to Problem 2 in~\cite{IdziakKrzaczkowski}.

\begin{thm} \label{theorem:tractability}
Assume that Conjecture~\ref{conj:CCcircuits} is true. Then, for every finite nilpotent algebra $\algA$ from a congruence modular variety $\CSAT(\algA)$ and $\CEQV(\algA)$ can be solved in quasipolynomial time $\mathcal O(2^{(\log n)^t})$ (where $t>0$ depends on $\algA$).
\end{thm}

\begin{proof}
We start with the equivalence problem $\CEQV(\algA)$. Without loss of generality we assume that $\algA$ contains a group operation $+$, such that $(A,+,0,-) = \prod_{i=1}^s \Z_{p_i}$ (otherwise we reduce to such an algebra by Theorem~\ref{theorem:coordinatisation} and Lemma~\ref{lemma:polynomialembedding}). When solving $\CEQV(\algA)$ it is sufficient to find an algorithm to check whether some input circuit is equivalent to the constant $0$-circuit (as $C = C'$ if and only if $C-C' = 0$). Thus we are only considering inputs $C$ and $0$ to $\CEQV(\algA)$.

By Proposition~\ref{proposition:main} we can identify $C(x_1,\ldots,x_n)$ with a $\CC^+[m]$-circuit $C'(y_1,\ldots,y_{sn})$  having a depth $d$ and $s$ output gates, where $d$ and $s$ only depend on $\algA$. Then $C(x_1,\ldots,x_n)$ is equivalent to $0$ if and only if $C'(c_1y_1,\ldots,c_s y_{sn})$ is constant and equivalent to $(0,0,\ldots,0)$ (see also the discussion after Proposition~\ref{proposition:main}). Let us denote the output gates of this circuit by $C_1,\ldots,C_s$.

In the case where one of the output gates $C_i$ does not compute the constant $0$ function, there is a tuple $\bar a = (a_1,\ldots,a_{sn})$ such $C_i(\bar a) \neq 0$. Let us pick $\bar a$ such that the number of coordinates $j$ with $a_j = 0$ is maximal. For the sake of simplicity, let us assume that $a_j = 0$ iff $j > k$ for some $k$. Then, the circuit $C_i(x_1,\ldots,x_k,0,\ldots,0)$ defines a nonconstant $0$-absorbing operation. Since we assume that Conjecture~\ref{conj:CCcircuits} is true, we have that the size $|C_i|$ of this circuit is bigger than $2^{k^q}$ for some $q > 0$. In other words $k \leq \log (|C_i|)^{l}$, with $l = q^{-1}$.

Thus in order to check, whether $C(\bar x)$ is equivalent to $0$, we only need to check, whether all $C_i$ evaluate to $0$ for all tuples that have at most $\log (|C_i|)^{l}$-many non-0 elements. There are $\binom{|C_i|}{\log (|C_i|)^{l}} = \mathcal O(|C_i|^{\log (|C_i|)^{l}})$ many such tuples. Since $|C_i|$ is linear in the size of $C$ we obtain an algorithm that runs in time $\mathcal O (|C|^{\log(|C|)^l}) = \mathcal O (2^{\log(|C|)^{l+1}})$. Thus the lemma holds for $t = l+1$.

The satisfiability problem $\CSAT(\algA)$ can similarly be reduced to checking whether the bounded depth $\CC^+[m]$-circuit $C'(c_1y_1,\ldots,c_s y_{sn})$ outputs the $s$-ary tuple $(0,0,\ldots,0)$ for some input $\bar a$. Let $f \colon \Z_m^s \to \Z_m$ be the function that outputs $0$ if and only if $x \neq (0,0,\ldots,0)$, and $1$ else. By Lemma~\ref{lemma:CC+fct} (2), $f$ can be computed by a $\CC^+[m]$-circuit of essential depth $3$. So clearly $C(\bar x) = 0$ has no solution, if and only if $f(C'(c_1y_1,\ldots,c_sy_{sm}))$ is constant and equal to $0$. Thus we reduced $\CSAT(\algA)$ to the equivalence problem for $\CC^+[m]$-circuit of a fixed depth, which has a quasipolynomial algorithm by the above. Note that the depth of this circuit is bigger than in the proof for $\CEQV$, thus we might obtain a bigger value for the constant $t$.
\end{proof}

We remark that already~\cite{barrington-eq-monoids} used the same approach to show that the program satisfiability problem over solvable monoids has a quasipolynomial algorithm (supposed that Conjecture~\ref{conj:CCcircuits} holds).  Moreover the algorithm described in Theorem~\ref{theorem:tractability} was already studied for the special case where $|A|$ is power of some prime $p$. As already mentioned, then $\CC[p]$-circuits of bounded depth cannot compute $\AND$, which corresponds to $\algA$ only having nonconstant $0$-absorbing polynomials of arity smaller than some constant $k$. So, in this case we only need to evaluate the circuits at tuples with $k$-many non $0$-entries, which gives us an algorithm that runs in polynomial time $\mathcal O(|C|^k)$. For $\CEQV$ this was observed in~\cite{AichingerMudrinski}, algorithms for $\CSAT$ were constructed in~\cite{kompatscher-supernilpotent}, \cite{IdziakKrzaczkowski} and~\cite{aichinger-eqsupernilpotent}.

Finite nilpotent algebras of prime power order and their direct products are provably the only finite nilpotent algebras, where we have a bound on the arity of nonconstant $0$-absorbing polynomials. Thus the algorithm described in Theorem~\ref{theorem:tractability} cannot be refined to run in polynomial time for general nilpotent algebras. However we remark that there are examples of 2-nilpotent, not supernilpotent algebras for which we can obtain other polynomial-time algorithms: It was shown in~\cite{IdziakKrzaczkowskiKawalek} that for nilpotent algebras $\algL \otimes \algU$, where both $\algL$ and $\algU$ are polynomially equivalent to finite vector space, $\CEQV$ and $\CSAT$ are in $\comP$. For $\CEQV$ this was generalized to all 2-nilpotent algebras in~\cite{KKK}.

At last we show that, under the assumption that there is an efficient way of \emph{uniformly} computing $\AND$ by bounded depth $\CC[m]$-circuits, there is a nilpotent algebra with intractable $\CSAT$ and $\CEQV$ problem (note however that this is a stronger assumption than $\AND$ being in $\CC^0$).

\begin{thm} \label{theorem:hardness}
Assume that there is a family $(C_n)_{n \in \mathbb N}$ of $\CC[m]$-circuits of depth $d$, that defines $\AND$, and is enumerable by a polynomial time Turing machine (i.e.  it  outputs $C_n(x_1,\ldots,x_n)$ on input $1^n$). Then there exists a finite nilpotent algebra $\algA$ such that $\CSAT(\algA)$ is $\comNP$-complete and $\CEQV(\algA)$ is $\comcoNP$-complete.
\end{thm}

\begin{proof}
Since $(C_n)_{n \in \mathbb N}$ defines $\AND$, $m$ cannot be a prime. In particular $m > 2$. We are first going to reduce the graph-colouring problem with $m$-colors to the circuit satisfiability problem for $\CC^+[m]$-circuits of depth bounded by $d+1$.

So let $G = (V,E)$ a graph. We are going to identify every vertex $v \in V$ with a variable $x_v$ over $\Z_m$, representing the color of $v$. Note that the $\CC^+[m]$-circuit $1 - \MOD_m (x_v - x_w)$ outputs $1$ if $x_v \neq x_w$, and $0$ else. Thus, if we define the circuit $H((x_v)_{v \in V}) = C_{|E|}((1 - \MOD_m (x_v - x_w))_{(v,w) \in E})$, then $H$ outputs $1$ if the assignment $v \mapsto x_v$ is a proper coloring of the graph and $0$ else. So $G$ is yes-instance of the $m$-coloring problem if and only if $H((x_v)_{v \in V}) = 1$ has a solution.

By our assumption $H$ can be computed in polynomial time from $G$, thus we reduced $m$-coloring to the satisfiability problem for $\CC^+[m]$-circuits of depth bounded by $d+1$. Furthermore, note that there is no $m$-coloring of $G$ if and only if $H((x_v)_{v \in V})$ is constant and equal to $0$. Therefore the complement of the $m$-coloring problem, reduces to checking the equivalence of $CC^+[m]$-circuits of depth $d+1$.

By Lemma~\ref{lemma:CCtonilpotent} there is a nilpotent algebra $\algA$, such that we can encode every $\CC^+[m]$-circuit $H$ of depth $d+1$ by a circuit $H'$ over $\algA$. Thus $m$-coloring reduces to $\CSAT(\algA)$ and its complement to $\CEQV(\algA)$.
\end{proof}

The proof of Theorem~\ref{theorem:hardness} uses a reduction of the coloring problems to $\CSAT$ with the help of $0$-absorbing polynomials; we remark that this idea was used in several other results in the relevant literature~\cite{BurrisLawrence-rings, goldmannrussell, HorvathSzabo-extendedgroups}.  

If $\AND$ satisfies the criterion in Theorem~\ref{theorem:hardness}, we  also say that $\AND$ is in \emph{$\comP$-uniform $\CC^0$}, for short. Note that, by Theorem~\ref{theorem:main}, this is equivalent to a statement about nilpotent algebras:

\begin{cor} \label{corollary:ANDinuCC0}
$\AND$ is in $\comP$-uniform $\CC^0$ if and only if there is a finite nilpotent algebra $\algA$ from a congruence modular variety, a family of nonconstant $0$-absorbing circuits $(C_n(x_1,\ldots,x_n))_{n \in \N}$ over $\algA$, and a polynomial time Turing machine that outputs $C_n(x_1,\ldots,x_n)$ on input $1^n$.
\end{cor}

\begin{proof}
This follows directly from the proof of Theorem~\ref{theorem:main}, since we showed how to compute an $n$-ary circuit computing $\AND$ from a circuit $(C_n(x_1,\ldots,x_n))_{n \in \N}$ over $\algA$ in linear time; and vice-versa.
\end{proof}

In Corollary~\ref{corollary:ANDinuCC0} it is essential that the nilpotent algebra is finite.  If we however also look at nilpotent algebras with infinite signature, it is quite easy to construct examples that have  nonconstant $0$-absorbing circuits, which are enumerable in polynomial time.  Theorem~6.1 in~\cite{IdziakKrzaczkowskiKawalek} gives an example of a 2-nilpotent algebra $\algA = (A,(f_n)_{n \in \mathbb N})$, for which this is the case,  and which has hard $\CSAT$ and $\CEQV$ problems. Interestingly this algebra $\algA$ is also polynomially equivalent to a finite nilpotent algebra $\algA'$, which in turn however fails to efficiently enumerate the same polynomial operations and satisfies $\CSAT(\algA') \in \comP$ and $\CEQV(\algA') \in \comP$.

\section*{Acknowledgments}

I would like to thank the anonymous reviewers for their many helpful remarks,  in particular for pointing out several inaccurate references.

\bibliographystyle{alphaurl}
\bibliography{notes_michael}

\end{document}